\documentclass[11pt]{article}
\usepackage{amsxtra,amssymb,amsthm,amsmath,latexsym}
\def\oH{\buildrel\circ\over H}
\def\oH1{\buildrel\circ\over H\kern-.02in{}^1}

\def\ep{\epsilon}

\begin{document}

\title{
Dynamical Systems Method for solving  equations with non-smooth monotone 
operators
\footnote{Math subject classification: 34R30,  35R25, 35R30,
37C35, 37L05, 37N30, 47A52,  
47J06, 65M30, 65N21$\quad$ key words: dynamical systems method, 
ill-posed problems, monotone operators, iterative methods, fixed-point 
maps }
}

\author{ A.G. Ramm\\
Mathematics Department, Kansas State University, \\
Manhattan, KS 66506-2602, USA\\
E:mail: ramm@math.ksu.edu
}

\date{}
\maketitle

\begin{abstract}
 Consider an operator equation (*) $B(u)+\ep u=0$ in a real Hilbert 
space,
where $\ep>0$ is a small constant.
The DSM (dynamical systems method) for solving equation (*)
consists of a construction of a Cauchy 
problem, which has the following properties:
 1) it has a global solution for an arbitrary initial data,\\
 2) this solution tends to a limit as time tends to infinity,\\
 3) the limit solves the equation $B(u)=0$.\\
Existence of the unique solution is proved by the DSM for 
equation
 (*) with monotone hemicontinuous operators $B$ defined on all of $H$.
If $\ep=0$ and equation (**) $B(u)=0$ is solvable, the DSM yields a 
solution to (**).
\end{abstract}

\section{Introduction}
In this paper a version of the DSM, dynamical systems method, is proposed
for
solving nonlinear operator
equation of the form:
$$
B(v)+\ep v=0,\quad \ep=const>0,
\eqno{(1)}
$$
where the operator $B: H\to H$ is a nonlinear monotone map in a Hilbert 
space $H$, and, in contrast to our earlier work [2], this operator may be
non-smooth.

The notions related to monotone operators, used in this paper, one finds,
e.g., in [1]. The DSM is applied to solving operator equations in [2],
where the nonlinear mapping $B$ was assumed locally twice  Fr\'echet
differentiable. In this paper existence of the Fr\'echet derivative of
$B(u)$ is not assumed. The novel feature in this paper is the
justification of the DSM for non-smooth monotone operators, and the novel
technique consists of the estimation of the derivative of the solutions 
to equations (2) and (7).

We make the following assumptions: 

A)\,\,{\it  
$B$ is a monotone, possibly nonlinear, hemicontinuous, defined on all of 
$H$ operator in a real Hilbert space $H$.}

If A) holds, then the set 
 $N:=\{z: B(z)=0\}$, if it is  non-empty, is closed and convex, and 
therefore it has the unique element $y$  with minimal norm.

Let $\dot u$ denote the derivative with respect to time. Consider the 
dynamical system ( that is, the Cauchy problem ):
$$
\dot w=-B(w)-\ep w, \,\,\, w(0)=w_0,
\eqno{(2)}
$$
where $w_0$ is arbitrary.

{\it The DSM  in this paper consists of solving
equation (1) by  
solving (2), and proving that for any initial approximation $w_0$ the 
following results (3) and (4) hold:}
$$
\exists w(t) \forall t>0, \quad \exists V_\ep:=w(\infty):=\lim_{t\to 
\infty}w(t),
\quad B(V_\ep)+\ep V_\ep=0, 
\eqno{(3)}
$$
and
$$\lim_{\ep\to 0}||V_\ep-y||=0.
\eqno{(4)}
$$
Conclusion (4) is known, but we give in subsection 2.5 a simple
proof for convenience of the reader.

It is also known that equation (1) under the assumtions A) has
a solution and this solution is unique. We prove this known fact by the
new method, the DSM. If $\ep=0$ then the limiting equation
$$
B(u)=0
\eqno{(5)}
$$
may have no solution (e.g., $B(u)=e^u$).We prove that if (5) has a
solution, then the DSM allows one to construct a solution to (5).
We assume 
$$
\ep (t)=\frac {c_1}{(c_0+t)^b}, \quad c_0>0,\,\,c_1>0,\,\, 0<b<1,
\eqno{(6)}
$$
where $c_0,c_1$ and $b$ are constants.
Consider the problem:
$$ 
\dot u=-B(u)-\ep(t) u, \,\,\, u(0)=u_0.
\eqno{(7)}
$$ 
Our results are stated in two theorems:

{\bf Theorem 1.} {\it If assumptions A) hold, equation (5) is solvable, 
and (6) holds, then 
problem (7) has a unique global solution $u(t)$, there exists strong limit
$u(\infty):=\lim_{t\to \infty}u(t)$, $B(u(\infty))=0$, and
$u(\infty):=y$ is the unique minimal-norm element in the set of all
solutions to (5).}

{\bf Theorem 2.} {\it If assumptions A) hold and $\ep=const>0$, then 
problem (2) has a unique global solution $w(t)$, there exists strong limit 
$w(\infty)$, and $w(\infty)$ solves (1).}

In Section 2 proofs are given.

\medskip

 \section{Proofs}

1. In Lemma 1 the existence of the unique global solutions
to problems (2) and (7) is claimed.
The proof of this Lemma will be given at the end of the paper
to make the presentation self-contained. The result of Lemma 1 is
known. Our proof is based on the Peano approximations (cf [1]).

{\bf Lemma 1.} {\it If $\ep=const>0$ and assumptions A)  hold, then
problem (2) has a unique global solution. If assumptions A) and (6) hold,
then  (7) has a unique global solution.}
  
2. {\bf Proof of Theorem 2.} The proof consists of the following steps:

a) we prove: 
$$
\sup_{t\geq 0}||u(t)||<c<\infty; \quad g(t)\leq g(0)e^{-\ep t}, \quad 
g(t):=||w(t+h)-w(t)||,
\eqno{(8)}
$$
where $c>0$ stands for various estimation constants, and $h>0$
is an arbitrary number.

 This and the Cauchy test imply the existence of $V_\ep:=w(\infty)$.

b) we prove that 
$$
||\dot w(t)||\leq ||\dot w(0)||e^{-\ep t}.
\eqno{(9)}
$$
 Thus, $\lim_{t\to \infty}||\dot w(t)||=0$. Estimate (9) implies
$\int_0^\infty ||\dot w(t)||dt<\infty$. This again, independently, implies
the existence of $V_\ep:=w(\infty)$.

From a), b), and from (2), one concludes that $V_\ep$ solves equation (1).
To pass to the limit in (2) one uses demicontinuity of hemicontinuous,
monotone, defined on all of $H$, operators, i.e., the property which
says that $w\to V$
implies $B(w)\rightharpoonup B(V)$, (cf,  e.g., [1], p.98).
Here and below $\rightharpoonup$ denotes weak convergence in $H$.

The proof of (4) is given in subsection 2.5.

To complete the proof of Theorem 2, let us prove (8) and (9).

Let $z:=w(t+h)-w(t)$ and $g:=||z||.$
From equation (2) and from the monotonicity of $B$ one gets: 
$$
g \dot g= -(B(w(t+h))-B(w(t))+\ep z, z)\leq -\ep g^2.
$$
Since $g\geq 0$, this implies the second half of (8).
Its first half, namely the estimate  
$\sup_{t\geq 0}||u(t)||<c<\infty$, is proved below formula (16)
for the more general case when $\epsilon$ depends on $t$.

Let $\psi:=\frac {||z||}h.$ Then, as above, one gets 
$\psi \dot \psi\leq -\ep \psi^2.$
Thus, $\psi (t)\leq \psi(0)e^{-\ep t}.$ 
Let $h\to 0$ and get (9). Theorem 2 is proved. $\Box$

$$ $$ 
3. {\bf Proof of Theorem 1.} 
The scheme of the proof is similar to the one used above, but there are 
new points due to the dependence of $\ep(t)$ on $t$ now.
 Denote $g(t):=||u(t+h)-u(t)||$ and  $z:=u(t+h)-u(t)$. From (7)
one gets
$$
g \dot g= -(B(w(t+h))-B(w(t))+\ep(t) z,
z)-(\ep(t+h)-\ep(t))(u(t+h),u(t))\leq 
$$ 
$$
-\ep g^2 +|\ep(t+h)-\ep(t)|||u(t+h)||g.
\eqno{(10)}
$$
We prove below that
$$
\sup_{t>0}||u(t)||\leq c<\infty, \quad c=const>0.
\eqno{(11)}
$$
 From (10), (11) and (6) one gets the following  differential inequality:
$$
\dot g\leq -\ep (t)g+hc|\dot \ep(t)|,
\eqno{(12)}
$$
where $c$ is defined in (11).

From (12) one gets:
$$
g(t)\leq e^{-\int_0^t\ep(s)ds}[g(0)+hc \int_0^te^{\int_0^s\ep(x)dx}|\dot
\ep(s)|ds].
\eqno{(13)}
$$
From (6) and (13) one gets 
$$\lim_{t\to \infty}g(t)=0, \quad \forall h>0. 
\eqno{(14)}
$$
Indeed, if $a(t):=e^{\int_0^t\ep(s)ds},$ then 
$a^{-1}(t)\int_0^t a(s)|\dot \ep (s)|ds=O(\frac 1 t)$ as $t\to \infty$,
as one derives from assumption (6).

From (11) it follows that there exists a sequence $t_n\to \infty$,
such that $u(t_n)\rightharpoonup v$, where $v\in H$ is some element.
We prove below that $B(v)=0$ by passing to the limit
$t_n\to \infty$ in equation (7), using assumption (6), inequality 
(11), and relation (15), which one obtains 
dividing (14) by $h$ and letting $h\to 0$:
$$\lim_{t\to \infty}||\dot u(t)||=0.
\eqno{(15)}
$$
Passing to the limit $t_n\to \infty$ in (7), proves that $u(\infty):=v$
solves (5).

Let us prove (11). In this
proof we use the assumption that (5) has a solution $y$. 

Denote
$u(t)-y:=p(t)$ and $||p||:=q$. Then 
$$
\dot p=-(B(u)-B(y))-\ep(t)p-\ep(t)y.
\eqno{(16)}
$$
Multiplying this by $p$ and using the monotonicity of $B$, one gets:
$$
\dot q\leq -\ep (t)q +\ep(t) ||y||.
\eqno{(17)}
$$
This implies $q(t)\leq a^{-1}(t)[q(0)+||y||\int_0^t a(s) 
\epsilon(s)ds]$. Thus,
$$
||u(t)-y||:=q(t)\leq c,
\eqno{(18)}
$$
so (11) follows (with a different $c$).

Let us now prove the existence of the strong limit $u(\infty)$ and
the relation $u(\infty)=y$, where $y$ is the unique minimal-norm solution 
to (5).

From (11) it follows that there is a sequence $t_n\to \infty$ such that
$u(t_n)\rightharpoonup v$. From (6), (11), (15), and (7) one gets
$\lim_{n\to \infty}B(u(t_n))=0$. This and assumptions A) imply $B(v)=0$.

Let us prove that $u(t_n)\to v$. Since $u(t_n)\rightharpoonup v$, one
gets 
$$||v||\leq  \liminf_{n\to \infty}||u(t_n)||.$$ 
If
$\limsup_{n\to \infty}||u(t_n)||\leq ||v||$, then $\lim_{n\to 
\infty}||u(t_n)||=||v||$, and together with the weak convergence
$u(t_n)\rightharpoonup v$ this implies strong convergence 
$u(t_n)\to v$. 

To prove that $\limsup_{n\to \infty}||u(t_n)||\leq ||v||$,
we need some preparations. First, (6) implies that $\int_0^t\epsilon(s)ds
\sim \frac {t^a}a$ as $t\to \infty$, where $a:=1-b$, $0<a<1$. Second,
(13) implies $||\dot{u}(t)||\leq c/t$ as $t\to \infty$, where $c>0$ is a 
constant.
Indeed, if (6) holds, then $a^{-1}(t)\int_0^t a(s)|\dot \ep (s)|ds=O(\frac 
1 t)$ as
$t\to \infty$. Equations $B(v)=0$ and (7) imply
$(B(u(t_n))-B(v),u(t_n)-v) +\epsilon (t_n)(u(t_n), u(t_n)-v)=
-(\dot{u}(t_n), u(t_n)-v)$. Since $B$ is monotone, it follows that
$(u(t_n), u(t_n)-v)\leq \frac c{t_n \epsilon (t_n)}$.
Thus, $\limsup_{n\to \infty}||u(t_n)||\leq ||v||$, because
$\lim_{n\to \infty}t_n \epsilon (t_n)=\infty$.  

Let us prove that $v=y$, where $y$ is the unique minimal-norm solution to 
(5). Replacing $v$ by $y$ in the above argument yields
$(u(t_n), u(t_n)-y)\leq \frac c{t_n \epsilon (t_n)}$, so
$||v||=\limsup_{n\to \infty}||u(t_n)||\leq ||y||$. Since $y$ is the 
unique minimal-norm solution to (5), and $v$ solves (5), it follows that 
$v=y$.

Since the limit $\lim_{n\to \infty}u(t_n)=v=y$ is the same for every
subsequence $t_n\to \infty$, for which the weak limit
of $u(t_n)$ exists, one concludes that the strong limit
$\lim_{t\to \infty}u(t)=y$. Indeed, assuming that for some sequence 
$t_n\to \infty$ the limit of $u(t_n)$ does not exist, one selects
a subsequence, denoted again $t_n$, for which the weak limit
of $u(t_n)$ does exist, and proves as before that this limit is $y$,
thus getting a contradiction.
 Theorem 1 is proved. $\Box$

4. {\bf Proof of Lemma 1.} 

Let $F$ be a nonlinear map satisfying assumptions A). In Lemma 1 this map
is $F(u):=B(u)+\ep u$ for equation (2) and $F(u):=B(u)+\ep(t) u$ for
equation (7).  Our argument holds for any $F$ satisfying assumptions A).
We want to prove that the problem
$$
\dot w=-F(w), \quad w(0)=w_0,
\eqno{(19)}
$$
has a unique global solution.

Uniqueness of the solution is immediate: if $w$ and $v$ are solutions to
(19), and $z:=w-v$, then $\dot z= -[F(w)-F(v)],$ $z(0)=0$. Multiplying by
$z$ and using the monotonicity of $F$, one gets $(\dot z,z)\leq 0$,
so $||z(t)||\leq 0$, and the uniqueness follows.

The proof of the global existence is less simple. To make
the paper selfcontained let us give a simplified version of the proof 
(cf [1]).

Consider the equation:
$$
w_n(t)=w_0-\int_0^t F(w_n(s-\frac 1n))ds,
\,\,\, t>0; \quad w_n(t)=w_0, \,\,\, t\leq 0.
\eqno{(20)}
$$
We wish to prove that 
$$
\lim_{n\to \infty}w_n(t)=w(t),\quad \forall t>0,
\eqno{(21)}
$$
 where $w$ solves (19).
Recall that assumptions A) imply demicontinuity of
$F$.

Fix an arbitrary $T>0$, and let $B(w_0,r)$ be the ball centered at $w_0$
with radius $r>0$. Let $\sup_{u\in B(w_0,r)}||F(u)||:=c.$
Then (20) implies $||w_n(t)-w_0||\leq ct$. If $t\leq r/c$, then $w_n(t)\in 
B(w_0,r)$, and $||\dot w_n(t)||\leq c.$ Define 
$$z_{nm}(t):=w_n(t)-w_m(t),\quad ||z_{nm}(t)||:=g_{nm}(t).
$$
From (20) one gets:
$$
g_{nm}\dot g_{nm}= -(F(w_n(t-\frac 1n))-F(w_m(t-\frac 1m)),
w_n(t)-w_m(t)):=I. 
$$
One has:
$$
I=-(F(w_n(t-\frac 1n))-F(w_m(t-\frac 1m)),w_n(t-\frac 1n)-w_m(t-\frac 1m))
$$
$$
-(F(w_n(t-\frac 1n))-F(w_m(t-\frac 1m)),w_n(t)-w_n(t-\frac
1n)-(w_m(t)-w_m(t-\frac 1m))).
$$
Using the monotonicity of $F$, the estimate $\sup_{w\in
B(w_0,r)}||F(w)||\leq c$,
and the estimate $||\dot w_n(t)||\leq c$, one gets: 
$$
I\leq 4c^2(\frac 1n+\frac 1m). 
$$
Therefore 
$$
 g_{nm}\dot g_{nm}\leq 4c^2(\frac 1n+\frac 1m)\to 0 \quad as \quad n,m\to 
\infty.
\eqno{(22)}
$$
This implies
$$
\lim_{n,m\to \infty}g_{nm}(t)=0, \quad 0\leq t\leq \frac r c.
\eqno{(23)}
$$
Therefore there exists the strong limit $w(t)$:
$$
\lim_{n\to \infty}w_{n}(t)=w(t), \quad 0\leq t\leq \frac r c.
\eqno{(24)}
$$
The function $w$, defined in (24), satisfies the integral equation:
$$
w(t)=w_0-\int_0^tF(w(s))ds,
\eqno{(25)}
$$
and solves problem (19). If $F$ is continuous, then problem (19) and
equation (25) are equivalent. If $F$ is demicontinuous, then they are also
equivalent, but the derivative in (19) should be understood in the weak
sense. We have proved the existence of the unique local solution to
(19).

To prove that the solution to (19) exists for any $t\in [0, \infty)$,
let us assume that the solution exists on $[0,T),$ but not on a larger
interval
$[0,T+d)$, $d>0$, and show that this leads to a contradiction. It is
sufficient
to prove that the finite limit:
$$\lim_{t\to T}w(t)
\eqno{(26)}
$$
does exist, because then one can solve locally, 
on the interval $[T, T+d)$, equation (19) with the
initial data $w(T)=\lim_{t\to T} w(t)$, and construct the 
solution to (19) on the interval $[0,T+d)$, thus getting a contradiction.

To  prove that the finite limit (26) exists, consider
$$
w(t+h)-w(t):=z(t), \quad ||z||:=g.
$$
One has $\dot z=-[F(w(t+h))-F(w(t))]$. Using the monotonicity of $F$, one
gets $(z,\dot z)\leq 0$. Thus, 
$$
||w(t+h)-w(t)||\leq ||w(h)-w(0)||.
\eqno{(27)}
$$ 
The right-hand side in (27) tends to zero as $h\to 0$. This, and the
Cauchy test imply the existence of the finite limit (26).

Lemma 1 is proved. $\Box$ 

5. {\bf Proof of (4).}

This proof requires the following lemmas in which assumptions A)
hold and are not repeated:

{\bf Lemma 2.} {\it If $y$ solves (5) and $V_\ep$ solves (3), then:
$$
||V_\ep||\leq ||y||.
\eqno{(28)}
$$
}
{\bf Lemma 3.} {\it If $v_n\rightharpoonup v$ and $B(v_n)\to f$,
then $B(v)=f$.
}

{\bf Lemma 4.} {\it If $v_n\rightharpoonup v$ and $||v_n||\leq ||v||$,
then $v_n \to v$.
}

Assuming these lemmas, let us prove (4). From (28) one gets 
$V_\ep\rightharpoonup v$ (by a subsequence denoted again $V_\ep$).
Equation (3) implies $B(V_\ep)\to 0$. Thus, Lemma 3 yields
$B(v)=0$. From $V_\ep\rightharpoonup v$ and from (26) one gets
$$||v||\leq \liminf_{\ep \to 0}||V_\ep|| \leq \liminf_{\ep \to 0}||V_\ep||
\leq ||y||.
$$
Therefore $v=y$, since the solution to the equation (5),
which has minimal norm, is unique, if A) holds.
The weak convergence $V_\ep\rightharpoonup y$,
inequality (28), and Lemma 4
imply (4).

Let us prove Lemmas 2-4.

{\bf Proof of Lemma 2.} One has $B(V_\ep)+\ep V_\ep-B(y)=0$. 
Multiply this by $V_\ep-y$ and use the monotonicity of 
$B$ to get $\ep (V_\ep, V_\ep-y)\leq 0$. Since $\ep>0$, inequality (28)
follows. $\Box$

{\bf Proof of Lemma 3.} The monotonicity of $B$ implies
$$
(B(v_n)-B(v-tz), v_n-v+tz)\geq 0 \quad\forall z\in H\,\,\forall t>0.
$$
Letting $n\to \infty$ one gets $(f-B(v-tz),tz)\geq 0$, so 
$$
(f-B(v-tz),z)\geq 0.
$$
Letting $t\to 0$, one gets $(f-B(v),z)\geq 0$ $\forall z$.
Thus, $B(v)=f$. Lemma 3 is proved. $\Box$

{\bf Proof of Lemma 4.} One has 
$
||v||\leq \liminf_{n\to \infty}||v_n||\leq \limsup_{n\to
\infty}||v_n||\leq ||v||.$
Thus, 
$\lim_{n\to \infty}||v_n||=||v||.$ 
This and the weak convergence $v_n\rightharpoonup v$, imply:
$||v-v_n||^2=||v_n||^2+||v||^2-2\Re (v_n,v)\to 0.$
Lemma 4 is proved. $\Box$.

Proof of (4) is completed. $\Box$

\end{document}